\documentclass[12pt,reqno]{amsart}

\usepackage{amssymb}

\numberwithin{equation}{section}

\newtheorem{theorem}{Theorem}[section]
\newtheorem{proposition}[theorem]{Proposition}
\newtheorem{lemma}[theorem]{Lemma}
\newtheorem{corollary}[theorem]{Corollary}

\theoremstyle{definition}
\newtheorem{example}[theorem]{Example}
\newtheorem{remark}[theorem]{Remark}

\begin{document}

\baselineskip=15pt

\title[Stability of pulled back parabolic vector bundles]{On the stability of pulled back 
parabolic vector bundles}

\author[I. Biswas]{Indranil Biswas}

\address{School of Mathematics, Tata Institute of Fundamental
Research, Homi Bhabha Road, Mumbai 400005, India}

\email{indranil@math.tifr.res.in}

\author[M. Kumar]{Manish Kumar}

\address{Statistics and Mathematics Unit, Indian Statistical Institute,
Bangalore 560059, India}

\email{manish@isibang.ac.in}

\author[A. J. Parameswaran]{A. J. Parameswaran}

\address{School of Mathematics, Tata Institute of Fundamental
Research, Homi Bhabha Road, Mumbai 400005, India}

\email{param@math.tifr.res.in}



\date{}

\begin{abstract}
Take an irreducible smooth projective curve $X$ defined over an algebraically closed 
field of characteristic zero, and fix finitely many distinct point $D\, =\, \{x_1,\, 
\cdots,\, x_n\}$ of it; for each point $x\, \in\, D$ fix a positive integer $N_x$. Take a 
nonconstant map $f\, :\, Y\, \longrightarrow \, X$ from an irreducible smooth projective 
curve. We construct a natural subbundle $\mathcal{F}\, \subset\, f_*{\mathcal O}_Y$ using 
$(D,\, \{N_x\}_{x\in D})$. Let $E_*$ be a stable parabolic vector bundle whose parabolic 
weights at each $x\, \in\, D$ are integral multiples of $\frac{1}{N_x}$. We prove that 
the pullback $f^*E_*$ is also parabolic stable, if ${\rm rank}(\mathcal{F})\,=\, 1$. 
\end{abstract}
\maketitle

\section{Introduction}

We begin by recalling the main result of \cite{BP}. Let $f\, :\, Y\, \longrightarrow\,
X$ be a surjective separable morphism between irreducible smooth projective curves
defined over an algebraically closed field $k$. It is called genuinely ramified if
the rank of the maximal semistable subbundle $F\, \subset\, f_*{\mathcal O}_Y$ is one.
The main result of \cite{BP} says that the pullback $f^*E$ of every stable vector bundle
$E$ on $X$ is also stable, provided $f$ is genuinely ramified.

Our aim here is to prove an analogue of it for parabolic vector bundles, but under
an extra assumption that the characteristic of the base field $k$ is zero.

Let $X$ be an irreducible smooth projective curve defined over an algebraically closed 
field $k$ of characteristic zero. Fix finitely many points $D\,=\, \{x_1,\, \cdots,\, 
x_n\}\,\subset\, X$, and for each $x\, \in\, D$ fix a positive integer $N_x$. We consider 
the category of parabolic vector bundles $E_*$ on $X$ with parabolic divisor $D$ such 
that all the parabolic weights of $E_*$ at any $x\, \in\, D$ are integral multiples of 
$\frac{1}{N_x}$.

Let $f\, :\, Y\, \longrightarrow\, X$ be a nonconstant morphism from an irreducible 
smooth projective curve $Y$. Using the above data $(D,\, \{N_x\}_{x\in D})$ we
construct a natural subbundle
\begin{equation}\label{ui}
\mathcal{F}\, \subset\, f_*{\mathcal O}_Y\,.
\end{equation}
This $\mathcal{F}$ is an analogue 
of the above mentioned maximal semistable subbundle $F\, \subset\, f_*{\mathcal O}_Y$ in the
context of parabolic bundles. It should be clarified that $f_*{\mathcal O}_Y$ is
parabolic polystable (see Proposition \ref{prop1}), so $\mathcal{F}$ is not related
to the Harder--Narasimhan filtration or the socle filtration of $f_*{\mathcal O}_Y$.

We prove the following (see Theorem \ref{thm3}):

\begin{theorem}\label{thmi}
Let
$$
f\, :\, Y\, \longrightarrow\, X
$$
be a nonconstant map between irreducible smooth projective curves defined over an 
algebraically closed field $k$
of characteristic zero. Take any stable parabolic vector bundle $E_*$ on $X$ with 
parabolic structure over $D$ such that all the parabolic weights of $E_*$ at each point 
$x\, \in\, D$ are integral multiples of $\frac{1}{N_x}$. If the rank of the vector bundle 
${\mathcal F}\,\longrightarrow\, X$ in \eqref{ui} is one, then the parabolic
vector bundle $f^*E_*$ on $Y$ is also stable.
\end{theorem}

We also prove the following converse of Theorem \ref{thmi} (see Lemma \ref{lem4}):

\begin{lemma}\label{lem-i}
Let
$$
f\, :\, Y\, \longrightarrow\, X
$$
be a nonconstant map between irreducible smooth projective curves defined over an
algebraically closed field $k$ of characteristic zero.
Assume that the rank of the holomorphic vector bundle
${\mathcal F}\,\longrightarrow\, X$ in \eqref{ui} is at least two.
Then there is a stable parabolic vector bundle $E_*$ on $X$ with parabolic structure over $D$
such that
\begin{enumerate}
\item all the parabolic weights of $E_*$ at each $x\, \in\, D$ are integral multiples
of $\frac{1}{N_x}$, and

\item the parabolic vector bundle $f^*E_*$ on $Y$ is not stable.
\end{enumerate}
\end{lemma}

Theorem \ref{thmi} is first proved for $k\,=\, \mathbb C$; see Theorem \ref{thm2}.
Then theorem \ref{thmi} is deduced using Theorem \ref{thm2}.

\section{Direct image and parabolic structure}

Let $k$ be an algebraically closed field of characteristic zero.

Let $X$ be an irreducible smooth projective curve
defined over $k$. Fix a nonempty finite subset
$$
D\, :=\, \{x_1,\, \cdots,\, x_n\}\, \subset\, X
$$
The reduced effective divisor $x_1+\ldots + x_n$ on $X$ will also be denoted
by $D$. A quasiparabolic structure on a vector bundle $E$ on $X$
is a filtration of subspaces of the fiber $E_{x_i}$ of $E$ over $x_i$
$$
E_{x_i}\,=\,E_{i,1}\,\supset\, E_{i,2}\,\supset\, \cdots\,
\supset \,E_{i,l_i} \,\supset\, E_{i,l_i+1}\,=\, 0
$$
for every $1\, \leq\, i\, \leq\, n$. A parabolic structure on
$E$ is a quasiparabolic structure as above together with a string of rational numbers
$$
0\,\leq\, \alpha_{i,1} \,< \,\alpha_{i,2} \,<\,
\cdots \,<\, \alpha_{i,l_i}\,<\, 1
$$
for every $1\, \leq\, i\, \leq\, n$. The above number $\alpha_{i,j}$
is called the parabolic weight of the subspace $E_{i,j}$.
(See \cite{MS}, \cite{MY}, \cite{Bh}, \cite{IIS}, \cite{In}.)

A parabolic vector bundle is a vector bundle $E$ equipped with a parabolic structure
$(\{E_{i,j}\},\, \{\alpha_{i,j}\})$. For notational convenience,
$(E,\, (\{E_{i,j}\},\, \{\alpha_{i,j}\}))$ will be denoted by $E_*$.
The divisor $D$ is known as the parabolic divisor of $E_*$.

The parabolic degree of $E_*$ is defined to be
$$
\text{par-deg}(E_*)\, :=\, \text{degree}(E)\,+\,\sum_{i=1}^n\sum_{j=1}^{l_i} \alpha_{i,j}
\cdot\dim (E_{i,j}/E_{i,j+1})\, ,
$$
and $\frac{\text{par-deg}(E_*)}{\text{rank}(E)}\, \in\, \mathbb Q$
is denoted by $\text{par-}\mu(E_*)$.

Let $F\, \subset\, E$ be a subbundle. Then a parabolic structure on $E$ produces a 
parabolic structure on $F$. The parabolic divisor for the induced parabolic structure on 
$F$ is $D$ itself. A subspace $0\, \not=\,V\, \subset\, F_{x_i}$ appears in the 
quasiparabolic filtration of $F_{x_i}$ if $$F_{x_i}\bigcap E_{i,j+1}\, \not=\, V\,=\, 
F_{x_i}\bigcap E_{i,j}$$ for some $1\,\leq\, j\, \leq\, l_i$. The parabolic weight of 
such a subspace $V$ is $\alpha_{i,j}$. The vector bundle $F$ with the induced parabolic 
structure will be denoted by $F_*$.

The parabolic bundle $E_*$ is called \textit{stable} (respectively, \textit{semistable}) if
$$
\text{par-}\mu(F_*)\, <\, \text{par-}\mu(E_*) \ \ \text{ (respectively, }\
\text{par-}\mu(F_*)\, \leq\, \text{par-}\mu(E_*)\text{)}
$$
for every subbundle $0\, \not=\, F \, \subsetneq \, E$. The parabolic bundle $E_*$ is
called \textit{polystable} if the following two conditions hold:
\begin{itemize}
\item $E_*$ is parabolic semistable, and

\item $E_*$ is a direct sum of stable parabolic bundles.
\end{itemize}

Let $Y$ be an irreducible smooth projective curve and
\begin{equation}\label{e2}
f\, :\, Y\, \longrightarrow\, X
\end{equation}
a nonconstant morphism. Let
\begin{equation}\label{e1}
D_f\, :=\, \{p_1,\, \cdots,\, p_m\}\, \subset\, X
\end{equation}
be the subset over which the map $f$ is ramified.

\begin{proposition}\label{prop1}
The direct image $f_*{\mathcal O}_Y$ has a parabolic structure whose parabolic 
divisor is $D_f$ defined in \eqref{e1}. This parabolic bundle given by
$f_*{\mathcal O}_Y$ is parabolic polystable of parabolic degree zero.
\end{proposition}

\begin{proof}
We use a local model of the map $f$ to describe the parabolic structure on
$f_*{\mathcal O}_Y$. Take
\begin{itemize}
\item $Y$ to be an open subset of $U\, \subset \, {\mathbb A}^1_k$
containing $0$,

\item $X$ to be the image of $U$ under the map ${\mathbb A}^1_k\, \longrightarrow\,
{\mathbb A}^1_k$ defined by $z\, \longmapsto\, z^d$, where $d$ is a positive integer, and

\item $f$ to be the map $z\, \longmapsto\, z^d$.
\end{itemize}
Then the quasiparabolic filtration of the fiber $(f_*{\mathcal O}_U)_0\,=\,
(f_*{\mathcal O}_U)_{p_1}$ over $0\,=\, p_1$ is given by the
image of the fibers of the filtration of subsheaves
$$
f_*{\mathcal O}_U\,\supset\, f_*({\mathcal O}_U(-p_1))\,\supset\, f_*({\mathcal O}_U(-2p_1))\,
\supset\,\cdots\,\supset\,f_*({\mathcal O}_U(-(d-1)p_1))\,\supset\,f_*({\mathcal O}_U(-dp_1))\, ,
$$
and the parabolic weight of the image of the fiber $(f_*({\mathcal O}_U(-kp_1)))_{p_1}$
in $(f_*{\mathcal O}_U)_{p_1}$ is $\frac{k}{d}$. Note that the image of the fiber
$(f_*({\mathcal O}_U(-dp_1)))_{p_1}$ in $(f_*{\mathcal O}_U)_{p_1}$ is zero, because by the
projection formula we have $f_*({\mathcal O}_U(-dp_1))\,=\, (f_*{\mathcal O}_U)\otimes 
{\mathcal O}_{f(U)}(-p_1)$.

Now consider the map $f$ in \eqref{e2}. For each $z_i\, \in\, D_f$ (see \eqref{e1}),
let $\{y^1_i,\, \cdots,\, y^{b_i}_i\}$ be the reduced inverse image $f^{-1}(z_i)_{\rm red}$.
For $1\, \leq\, j\, \leq\, b_i$, let $U^j_i$ be the formal completion of $y^j_i$ in $Y$.
The restriction of $f$ to $U^j_i$ will be denoted by $f^j_i$.
Now we have
\begin{equation}\label{j1}
(f_*{\mathcal O}_Y)_{z_i}\, =\, \bigoplus_{j=1}^{b_i}
\left((f^j_i)_* {\mathcal O}_{U^j_i}\right)_{z_i}\, .
\end{equation}
Each direct summand $\left((f^j_i)_* {\mathcal O}_{U^j_i}\right)_{z_i}$ of the fiber 
$(f_*{\mathcal O}_Y)_{z_i}$ in \eqref{j1} has a parabolic structure which is described 
above. The parabolic structure on $(f_*{\mathcal O}_X)_{z_i}$ is given by the direct sum 
of the parabolic structures on the direct summands in \eqref{j1}.

We will give an alternative description of the parabolic structure on $f_*{\mathcal O}_Y$.
Let $Z$ be an irreducible smooth projective curve and
\begin{equation}\label{e3}
\phi\, :\, Z\, \longrightarrow\, Y
\end{equation}
a nonconstant morphism, such that $f\circ\phi\, :\, Z\, \longrightarrow\, X$ is a
(ramified) Galois covering. Let $\Gamma\,=\, \text{Gal}(f\circ\phi)\, \subset\,
\text{Aut}(Z)$ be the Galois group.

Let $k[\Gamma]$ denote the algebra of functions on the finite group $\Gamma$. The
left-translation action of $\Gamma$ on itself produces an action of $\Gamma$ on
$k[\Gamma]$. On the other hand, the group $\Gamma\, \subset\, \text{Aut}(Z)$ has
a tautological action on $Z$. Consider the diagonal action of $\Gamma$ on
$Z\times k[\Gamma]$. This action makes the trivial vector bundle
\begin{equation}\label{e4}
Z\times k[\Gamma]\, \longrightarrow\, Z
\end{equation}
a $\Gamma$--equivariant vector bundle on $Z$. Using the natural correspondence between 
equivariant bundles and parabolic bundles (see \cite{Bi1}, \cite{Bo1}, \cite{Bo2}), this 
$\Gamma$--equivariant vector bundle $Z\times k[\Gamma]$ on $Z$ produces a parabolic vector 
bundle on $Z/\Gamma\,=\, X$. This parabolic vector bundle on $X$ will be denoted by 
$W_*$. The parabolic divisor $\text{par-div}(W_*)$ for $W_*$ is the subset of $X$ over 
which the map $f\circ\phi$ is ramified, where $\phi$ is the map in \eqref{e3}. Note that 
$D_f\, \subset\, \text{par-div}(W_*)$, and $\text{par-div}(W_*)$ may be larger than 
$D_f$.

The vector bundle underlying the parabolic vector bundle $W_*$ is $(f\circ\phi)_*{\mathcal 
O}_Z$ \cite{Bi2}, \cite{Pa}. On the other hand, we have
$$
f_*{\mathcal O}_Y\, \subset\, (f\circ\phi)_*{\mathcal O}_Z\, .
$$
In fact, $f_*{\mathcal O}_Y$ is a subbundle of $(f\circ\phi)_*{\mathcal O}_Z$.

Hence the parabolic structure of $W_*$ induces a parabolic structure on $f_*{\mathcal O}_Y$.
Let
$$
(f_*{\mathcal O}_Y)_*\, \longrightarrow\, X
$$
denote the parabolic vector bundle with parabolic structure on $f_*{\mathcal O}_Y$
induced by $W_*$.

The parabolic structure on $(f_*{\mathcal O}_Y)_*$ over the complement 
$\text{par-div}(W_*)\setminus D_f$ is the trivial one, meaning $(f_*{\mathcal O}_Y)_*$ does 
not have any nonzero parabolic weight on the points of $\text{par-div}(W_*)\setminus D_f$.

It is straight-forward to check that $(f_*{\mathcal O}_Y)_*$ coincides with the parabolic 
bundle given by the parabolic structure on $f_*{\mathcal O}_Y$ constructed earlier. In
particular, the above parabolic structure on $(f_*{\mathcal O}_Y)_*$ does not
depend on the choice of the pair $(Z,\, \phi)$.

Since the vector bundle underlying the $\Gamma$--equivariant vector bundle in \eqref{e4} 
is polystable (it is in fact trivial), the corresponding parabolic bundle $(f_*{\mathcal 
O}_Y)_*$ is polystable \cite[p.~350--351, Theorem 4.3]{BBN}. Since the degree of the 
$\Gamma$--equivariant vector bundle in \eqref{e4} zero, it follows that the parabolic 
degree of $(f_*{\mathcal O}_Y)_*$ is zero \cite[p.~318, (3.12)]{Bi1}. \end{proof}

We refer the reader to \cite{Yo} for the definition of parabolic dual of parabolic
vector bundles.

\begin{lemma}\label{lem1}
The parabolic dual of the parabolic vector bundle $(f_*{\mathcal O}_Y)_*$,
constructed in the proof of Proposition \ref{prop1}, is $(f_*{\mathcal O}_Y)_*$ itself.
\end{lemma}

\begin{proof}
If $E_*$ is the parabolic vector bundle corresponding to an equivariant bundle
$V$, then the parabolic vector bundle corresponding to the equivariant bundle
$V^*$ is the parabolic dual $E^*_*$ of $E_*$ \cite{BBN}.

Let
\begin{equation}\label{e5}
k[\Gamma]\otimes k[\Gamma]\, \longrightarrow\, k
\end{equation}
be the pairing defined by $$\left(\sum_{\gamma\in \Gamma} a_\gamma \gamma,\,
\sum_{\gamma\in \Gamma} b_\gamma \gamma\right)\, \longmapsto\, \sum_{\gamma\in \Gamma}
a_\gamma b_\gamma\, ,$$ where $a_\gamma,\, b_\gamma\, \in\, k$.
Consider the $\Gamma$--equivariant vector bundle $Z\times k[\Gamma]\, \longrightarrow\, Z$
in \eqref{e4}. The paring in \eqref{e5} defines a homomorphism of coherent sheaves
$$
(Z\times k[\Gamma])\otimes (Z\times k[\Gamma])\, \longrightarrow\, {\mathcal O}_Z
$$
which is fiberwise nondegenerate. The resulting isomorphism of vector bundles
$$
Z\times k[\Gamma]\, \stackrel{\sim}{\longrightarrow}\, Z\times k[\Gamma]^*
\,=\,(Z\times k[\Gamma])^*
$$
is in fact $\Gamma$--equivariant. Therefore, we conclude that the parabolic
dual of $(f_*{\mathcal O}_Y)_*$ is $(f_*{\mathcal O}_Y)_*$ itself.
\end{proof}

\section{Construction of a parabolic subbundle}\label{se3}

Let $V_*\,=\, (V,\, (\{V_{i,j}\},\, \{\alpha_{i,j}\}))$ be a semistable parabolic bundle
on $X$ with parabolic divisor 
$$
\mathbb{D}\, :=\, \{t_1,\, \cdots,\, t_r\}\, \subset\, X\, .
$$
For any subbundle $F\, \subset\, V$, the parabolic vector bundle defined by $F$ equipped
with the parabolic structure induced by $V_*$ will be denoted by $F_*$.

For each parabolic point $t\, \in\,\mathbb{D}$, we fix an integer $N_t\, \geq\, 1$. Assume
that there is a subbundle
$$
F\, \subset\, V
$$
satisfying the following two conditions:
\begin{enumerate}
\item All the parabolic weights of $F_*$
at every $t\, \in\, \mathbb{D}$ are integral multiples of $\frac{1}{N_t}$. (If
$N_t\,=\, 1$, then $F_*$ does not have any nonzero parabolic weight at $t$.)

\item $\text{par-}\mu(F_*)\,=\, \text{par-}\mu(V_*)$.
\end{enumerate}

\begin{lemma}\label{lem2}
There is a unique maximal subbundle
$$
\mathcal{F}\, \subset\, V
$$
satisfying the following two conditions:
\begin{enumerate}
\item All the parabolic weights of $\mathcal{F}_*$
at every $t\, \in\, \mathbb{D}$ are integral multiples of $\frac{1}{N_t}$, and

\item $\text{{\rm par}-}\mu(\mathcal{F}_*)\,=\, \text{{\rm par}-}\mu(V_*)$.
\end{enumerate}
\end{lemma}

\begin{proof}
Let $F^1$ and $F^2$ be two subbundles of $V$ such that for $1\, \leq\, j\, \leq\, 2$,
\begin{enumerate}
\item all the parabolic weights of $F^j_*$
at every $t\, \in\, \mathbb{D}$ are integral multiples of $\frac{1}{N_t}$, and

\item $\text{par-}\mu(F^j_*)\,=\, \text{par-}\mu(V_*)$.
\end{enumerate}
Since $V_*$ is parabolic semistable, and $\text{par-}\mu(F^j_*)\,=\, \text{par-}\mu(V_*)$, it
follows immediately that $F^j_*$ is semistable for $j\,=\, 1,\, 2$. Consider the subsheaf
$F^1+ F^2\, \subset\, V$ equipped with the parabolic structure induced by the parabolic
structure of $V_*$; the resulting parabolic bundle will be denoted by $(F^1+F^2)_*$. So
\begin{equation}\label{z1}
\text{par-}\mu((F^1+F^2)_*)\,\leq \, \text{par-}\mu(V_*)
\end{equation}
because $V_*$ is parabolic semistable. On the other hand,
$(F^1+F^2)_*$ is a quotient of the direct sum $F^1_*\oplus F^2_*$, and
$F^1_*\oplus F^2_*$ is parabolic semistable with 
$$
\text{par-}\mu((F^1\oplus F^2)_*)\,= \, \text{par-}\mu(F^1_*)\,=\, \text{par-}\mu(F^2_*)\,=\,
\text{par-}\mu(V_*)\, .
$$
Hence we have
$$
\text{par-}\mu((F^1+F^2)_*)\,\geq \, \text{par-}\mu(V_*)\,.
$$
Combining this with \eqref{z1} we conclude that
\begin{equation}\label{z2}
\text{par-}\mu((F^1+F^2)_*)\,=\, \text{par-}\mu(V_*).
\end{equation}

We will show that $V/(F^1+ F^2)$ is torsionfree. To prove this, if $T_0$ is the torsion
part of $V/(F^1+ F^2)$, consider ${\mathbb S}\,=\, q^{-1}_0(T_0)$, where $q_0\,:=\,
V\, \longrightarrow\, V/(F^1+F^2)$ is the quotient map. Let ${\mathbb S}_*$ denote the
parabolic vector bundle given by $\mathbb S$ equipped with the parabolic structure induced
by the parabolic structure of $V_*$. If $T_0\, \not=\, 0$, then
$$
\text{par-}\mu({\mathbb S}_*)\,> \, \text{par-}\mu((F^1+F^2)_*)\,=\, \text{par-}\mu(V_*)
$$
(see \eqref{z2}). But this contradicts the given condition that $V_*$ is parabolic semistable.
Therefore, we conclude that $V/(F^1+ F^2)$ is torsionfree. In other words,
$$
F^3\, :=\, F^1+F^2\, \subset\, V
$$
is a subbundle.

Consider the parabolic vector bundle $F^3_*$ defined by $F^3$ equipped with the parabolic 
structure induced by the parabolic structure of $V_*$. Recall that for $1\, \leq\, j\, \leq\, 2$, all the parabolic 
weights of $F^j_*$ at every $t\, \in\, \mathbb{D}$ are integral multiples of 
$\frac{1}{N_t}$. This immediately implies that all the parabolic weights of $F^3_*$ at 
each $t\, \in\, \mathbb{D}$ are also integral multiples of $\frac{1}{N_t}$.

In \eqref{z2} we have seen that $\text{par-}\mu(F^3_*)\,=\, \text{par-}\mu(V_*)$.

Now take $\mathcal F$ to be the coherent subsheaf of $V$ generated by all subbundles
$$
F\, \subset\, V
$$
such that
\begin{enumerate}
\item all the parabolic weights of $F_*$ at each $t\, \in\,\mathbb{D}$ are integral
multiples of $\frac{1}{N_t}$, and

\item $\text{par-}\mu(F_*)\,=\, \text{par-}\mu(V_*)$.
\end{enumerate}
{}From the above observations on $F^3_*$ it follows immediately that this coherent
subsheaf ${\mathcal F}\, \subset\, V$ 
satisfies all the conditions in the statement of the lemma.
\end{proof}

As in \eqref{e2}, take any irreducible smooth projective curve $Y$ together with a
nonconstant morphism
$$
f\, :\, Y\, \longrightarrow\, X\, .
$$
As in \eqref{e1}, $D_f\, :=\, \{p_1,\, \cdots,\, p_m\}\, \subset\, X$
denotes the subset over which $f$ is ramified. Fix a divisor
\begin{equation}\label{a6}
D\, :=\, \{x_1,\, \cdots,\, x_n\}\, \subset\, X\, .
\end{equation}
Also, fix an integer
\begin{equation}\label{e6}
N_x\, \geq\, 1
\end{equation}
for each point $x\, \in\, D$.

\begin{proposition}\label{prop2}
Let $(f_*{\mathcal O}_Y)_*$ be the parabolic bundle defined by $f_*{\mathcal O}_Y$
equipped with the natural parabolic structure (see Proposition \ref{prop1}).
Then there is a unique maximal subbundle
$$
\mathcal{F}\, \subset\, f_*{\mathcal O}_Y
$$
satisfying the following three conditions:
\begin{enumerate}
\item For any $x\, \in\, D\bigcap D_f$, all the parabolic weights of $\mathcal{F}_*$ 
at $x$ are integral multiples of $\frac{1}{N_x}$ (see \eqref{e6}),

\item $\mathcal{F}_*$ does not have any nonzero parabolic weight over any
point of $D_f\setminus (D\bigcap D_f)$, and

\item $\text{{\rm par-deg}}(\mathcal{F}_*)\,=\, 0$.
\end{enumerate}
\end{proposition}

\begin{proof}
In Lemma \ref{lem2}, set
\begin{itemize}
\item $V_*\,=\, (f_*{\mathcal O}_Y)_*$, (so the parabolic divisor $\mathbb{D}$
in Lemma \ref{lem2} is now $D_f$),

\item $N_x\,=\, 1$ if $x\, \in\, D_f\setminus (D\bigcap D_f)$ (see \eqref{a6} for $D$), and

\item $N_x\,=\, N_x$ (see \eqref{e6}) if $x\, \in\, D\bigcap D_f$.
\end{itemize}
Recall from Proposition \ref{prop1} that $\text{par-deg}((f_*{\mathcal O}_Y)_*)\,=\, 0$.
So we have $$\text{par-}\mu((f_*{\mathcal O}_Y)_*)\,=\, 0.$$ Therefore, in view of Lemma
\ref{lem2} it suffices to show that there is a subbundle
$$
F\, \subset\, f_*{\mathcal O}_Y
$$
satisfying the following two conditions:
\begin{enumerate}
\item All the parabolic weights of $F_*$ at each $x\, \in\, D_f$ are integral multiples of
$\frac{1}{N_x}$, and

\item $\text{par-deg}(F_*)\,=\, 0$.
\end{enumerate}

Since
$$
H^0(Y,\, \text{Hom}({\mathcal O}_Y,\, {\mathcal O}_Y))\,=\,
H^0(Y,\, \text{Hom}(f^*{\mathcal O}_X,\, {\mathcal O}_Y))\,=\,
H^0(X,\, \text{Hom}({\mathcal O}_X,\, f_*{\mathcal O}_Y))
$$
(see \cite[p.~110]{Ha}), the identity map of ${\mathcal O}_Y$ produces a
nonzero homomorphism
$$
{\mathcal O}_X\, \hookrightarrow\, f_*{\mathcal O}_Y\, .
$$
This coherent subsheaf is actually a subbundle. Indeed, this follows immediately from the
fact that for any $\psi\, \in\, H^0(U,\, {\mathcal O}_U)$, where $U\, \subset\, X$
is a Zariski open subset, the section $\psi\circ f\, \in\,
H^0(f^{-1}(U),\, {\mathcal O}_{f^{-1}(U)})$ has the property that if $\psi (x)\, \not=\, 0$,
then $\psi\circ f$ does not vanish on any point of $f^{-1}(x)$.

{}From the construction of the parabolic structure on $f_*{\mathcal O}_Y$ in Proposition
\ref{prop1} it follows immediately that the induced parabolic weight on ${\mathcal O}_X$
at any $x\,\in\, D_f$ is zero. Consequently, $F\,=\, {\mathcal O}_X$ satisfies the above
two conditions. This proves the proposition.
\end{proof}

Equip the curve $X$ with the following orbifold structure: For each point $x\, \in\, D$ 
the inertia group is ${\mathbb Z}/N_x\mathbb Z$, where $N_x$ is the integer in 
\eqref{e6}. The curve $X$ equipped with this orbifold structure will be denoted by 
$\mathcal X$. An {\it \'etale covering}
\begin{equation}\label{e11}
\varphi\,\,:\,\, Z \longrightarrow\, \mathcal{X}
\end{equation}
is an irreducible smooth projective curve $Z$ together with a nonconstant morphism
\begin{equation}\label{e12}
\varphi_0\,\,:\,\, Z \longrightarrow\, X
\end{equation}
such that the following conditions hold:
\begin{itemize}
\item the map $\varphi_0$ is unramified over $X\setminus D\,=\,
X\setminus \{x_1,\, \cdots ,\, x_n\}$, and

\item for every $x\, \in\, D$, the order
of ramification of $\varphi_0$ at each $z\, \in\, \varphi^{-1}_0(x)$ is a divisor of $N_x$.
\end{itemize}
An \'etale covering $\varphi$ of $\mathcal X$ will be called \textit{nontrivial} if
${\rm degree}(\varphi_0)\, \geq\, 2$. 

\begin{theorem}\label{thm1}
Consider the map $f\, :\, Y\, \longrightarrow\, X$, and the
corresponding subbundle
$$
\mathcal{F}\, \subset\, f_*{\mathcal O}_Y
$$
in Proposition \ref{prop2}. Then the following two statements are equivalent:
\begin{enumerate}
\item There is a nontrivial \'etale covering
$$
\varphi\,\,:\,\, Z \longrightarrow\, \mathcal{X}
$$
(see \eqref{e11} and \eqref{e12})
and a morphism $\beta\,:\, Y\, \longrightarrow\, Z$, such that
$\varphi_0\circ\beta\,=\, f$.

\item The rank of $\mathcal F$ is bigger than one.
\end{enumerate}
\end{theorem}

\begin{proof}
First assume that there is a nontrivial \'etale covering
$$
\varphi\,\,:\,\, Z \longrightarrow\, \mathcal{X}\, ,
$$
and a morphism $\beta\,:\, Y\, \longrightarrow\, Z$, such that
$\varphi_0\circ\beta\,=\, f$. Consider the subbundle
\begin{equation}\label{te}
(\varphi_0)_*{\mathcal O}_Z\, \subset\, f_*{\mathcal O}_Y\, .
\end{equation}
The parabolic structure on $(\varphi_0)_*{\mathcal O}_Z$
constructed in Proposition \ref{prop1} coincides with the
one induced by the parabolic structure of $f_*{\mathcal O}_Y$
on the subbundle in \eqref{te}. The parabolic bundle defined by
this parabolic structure on $(\varphi_0)_*{\mathcal O}_Z$ will
be denoted by $((\varphi_0)_*{\mathcal O}_Z)_*$.

Using the given condition that $\varphi$ is an \'etale covering of $\mathcal X$
it is straight-forward to verify that for every point $x\, \in\, D$, all the
parabolic weights of $((\varphi_0)_*{\mathcal O}_Z)_*$ 
at $x$ are integral multiples of $\frac{1}{N_x}$. Since $\varphi_0$ is unramified over
the complement $X\setminus D$, the parabolic bundle $((\varphi_0)_*{\mathcal O}_Z)_*$
does not have any nonzero parabolic weights on $X\setminus D$.
Also, from Proposition \ref{prop1} we know that
$$
\text{par-deg}(((\varphi_0)_*{\mathcal O}_Z)_*)\,=\, 0.
$$
In view of these, from the uniqueness property of $\mathcal F$ in Proposition
\ref{prop2} we know that
\begin{equation}\label{e7}
(\varphi_0)_*{\mathcal O}_Z\, \subset\, {\mathcal F}\, .
\end{equation}
Since $\text{degree}(\varphi_0)\, \geq\, 2$, from \eqref{e7} we conclude that
$$
{\rm rank}(\mathcal{F})\, \geq\, 2\, .
$$

To prove the converse, assume that
\begin{equation}\label{e8}
{\rm rank}(\mathcal{F})\, \geq\, 2\, .
\end{equation}
As before, ${\mathcal F}_*$ denotes the parabolic bundle defined by $\mathcal F$ equipped with the
parabolic structure induced by the parabolic structure of $f_*{\mathcal O}_Y$.

The algebra structure of ${\mathcal O}_Y$ produces an algebra structure
\begin{equation}\label{e9}
\Phi\,\,:\,\, (f_*{\mathcal O}_Y)_*\otimes (f_*{\mathcal O}_Y)_*
\, \longrightarrow\, (f_*{\mathcal O}_Y)_*
\end{equation}
(see \cite{Yo}, \cite{BBN} for the tensor product of parabolic vector bundles).

Since the parabolic weights of $\mathcal{F}_*$ 
at every $x\, \in\, D$ are integral multiples of $\frac{1}{N_x}$,
we conclude that
the parabolic weights of $\mathcal{F}_*\otimes {\mathcal F}_*$ 
at $x\, \in\, D$ are also integral multiples of $\frac{1}{N_x}$. From the
given condition that $\text{par-deg}(\mathcal{F}_*)\,=\, 0$ it follows
immediately that
$$\text{par-deg}(\mathcal{F}_*\otimes {\mathcal F}_*)\,=\, 0\, .$$

We note that $\mathcal{F}_*\otimes \mathcal{F}_*$ is parabolic semistable
because $\mathcal{F}_*$ is so \cite[p.~346, Proposition 3.2]{BBN}.
Let $\Phi(\mathcal{F}_*\otimes \mathcal{F}_*)_*$ be the parabolic vector
bundle defined by $\Phi(\mathcal{F}_*\otimes \mathcal{F}_*)$ equipped with the
induced parabolic structure, where $\Phi$ is the homomorphism in \eqref{e9}. Since 
$\Phi(\mathcal{F}_*\otimes \mathcal{F}_*)_*\, \subset\, (f_*{\mathcal O}_Y)_*$
is a quotient parabolic bundle of the semistable
parabolic bundle $\mathcal{F}_*\otimes \mathcal{F}_*$, we have
$$
\text{par-}\mu(\Phi(\mathcal{F}_*\otimes \mathcal{F}_*)_*)\, \geq\,
\text{par-}\mu(\mathcal{F}_*\otimes \mathcal{F}_*)\,=\, 0\, .
$$
On the other hand,
$$
\text{par-}\mu(\Phi(\mathcal{F}_*\otimes \mathcal{F}_*)_*)\, \leq\,
\text{par-}\mu((f_*{\mathcal O}_Y)_*)\,=\, 0\, ,
$$
because $(f_*{\mathcal O}_Y)_*$ is polystable. Combining these, we have
$$
\text{par-}\mu(\Phi(\mathcal{F}_*\otimes \mathcal{F}_*)_*)\,=\,0.
$$
Since $\Phi(\mathcal{F}_*\otimes \mathcal{F}_*)_*$ is a quotient of
$\mathcal{F}_*\otimes \mathcal{F}_*$, and all the parabolic weights
of $\mathcal{F}_*\otimes \mathcal{F}_*$ at every $x\, \in\, D$ are integral
multiples of $\frac{1}{N_x}$, it follows that all the parabolic weights
of $\Phi(\mathcal{F}_*\otimes \mathcal{F}_*)_*$ at every $x\, \in\, D$ are also
integral multiples of $\frac{1}{N_x}$. 

Therefore, from the uniqueness property of $\mathcal F$ we conclude that
\begin{equation}\label{e10}
\Phi(\mathcal{F}_*\otimes \mathcal{F}_*)_*\, \subset\, \mathcal{F}_*\, .
\end{equation}

{}From \eqref{e10} it follows that there is a unique \'etale covering
\begin{equation}\label{u1}
\varphi\,\,:\,\, Z \longrightarrow\, \mathcal{X}\, ,
\end{equation}
where $\mathcal{X}$ is the orbifold in \eqref{e11}, and
a morphism $\beta\,:\, Y\, \longrightarrow\, Z$, such that
following two hold:
\begin{enumerate}
\item $\varphi_0\circ\beta\,=\, f$, and

\item the two subsheaves $(\varphi_0)_*{\mathcal O}_Z$ and $\mathcal{F}$
of $f_*\mathcal{O}_Y$ coincide.
\end{enumerate}
{}From \eqref{e8} and the above statement (2) we know that the \'etale covering
$\varphi$ in \eqref{u1} is nontrivial. This completes the proof.
\end{proof}

Consider the set-up of Theorem \ref{thm1}. Let
$$
Y'\, :=\, Y\setminus f^{-1}(D)\, \subset\, Y
$$
be the complement. Let
\begin{equation}\label{e13}
f'\, :=\, f\big\vert_{Y'}\, :\, Y'\, \longrightarrow\, \mathcal{X}
\end{equation}
be the restriction of $f$ to $Y'$. The \'etale fundamental groups of
$Y'$ and $\mathcal X$ will be denoted by $\pi_1(Y')$ and $\pi_1(\mathcal{X})$
respectively. The following two statements are evidently equivalent:
\begin{enumerate}
\item There is a nontrivial \'etale covering
$$
\varphi\,\,:\,\, Z \longrightarrow\, \mathcal{X}
$$
(see \eqref{e11} and \eqref{e12}) and a morphism $\beta\,:\, Y\,
\longrightarrow\, Z$, such that $\varphi_0\circ\beta\,=\, f$.

\item The homomorphism of \'etale fundamental groups
\begin{equation}\label{re3}
(f'_{\rm et})_*\, :\, \pi_1(Y')\, \longrightarrow\, \pi_1(\mathcal{X})
\end{equation}
induced by $f'$ in \eqref{e13} is not surjective.
\end{enumerate}

Therefore, Theorem \ref{thm1} gives the following:

\begin{corollary}\label{cor1}
Consider the map $f\, :\, Y\, \longrightarrow\, X$, and the
corresponding subbundle
$$
\mathcal{F}\, \subset\, f_*{\mathcal O}_Y
$$
in Proposition \ref{prop2}. Then the following two statements are equivalent:
\begin{enumerate}
\item The homomorphism of \'etale fundamental groups
$(f'_{\rm et})_*$ in \eqref{re3} is not surjective.

\item The rank of $\mathcal F$ is bigger than one.
\end{enumerate}
\end{corollary}

\section{Complex curves and pullback of stable parabolic bundles}

Throughout this section we assume that $k\,=\, \mathbb C$. The topological
fundamental group of any complex manifold or orbifold $\textbf{N}$ will be
denoted by $\pi^t_1(\textbf{N})$; this is to distinguish it from the \'etale
fundamental group of $\mathbf{N}$.

\subsection{Homomorphism of topological fundamental groups}

As before, $f\, :\, Y\, \longrightarrow\, X$ is a nonconstant holomorphic map
between irreducible complex projective curves; the map $f$ is ramified exactly over
$$D_f\,:=\, \{p_1,\, \cdots, \,p_m\}\,\subset\, X\, .$$ Fix an integer $N_x\, \geq\, 1$
for each $x\,\in\, D$, and the resulting orbifold is denoted by $\mathcal{X}$.
The curve $Y'$ and the map $f'$ are both as in \eqref{e13}.

\begin{proposition}\label{prop3}
The following two statements are equivalent:
\begin{enumerate}
\item the homomorphism of topological fundamental groups
\begin{equation}\label{r3}
f'_* \, :\, \pi^t_1(Y')\, \longrightarrow\, \pi^t_1(\mathcal{X})
\end{equation}
induced by $f'$ in \eqref{e13} is surjective.

\item The rank of the holomorphic vector bundle ${\mathcal F}\,\longrightarrow\, X$ 
in Proposition \ref{prop2} is one.
\end{enumerate}
\end{proposition}

\begin{proof}
First assume that the homomorphism
$$
f'_*\, :\, \pi^t_1(Y')\, \longrightarrow\, \pi^t_1(\mathcal{X})
$$
is surjective. The group $\pi_1(Y')$ (respectively, $\pi_1(\mathcal{X})$)
is the profinite completion of $\pi^t_1(Y')$ (respectively, $\pi^t_1(\mathcal{X})$).
Therefore, from the surjectivity of the above homomorphism $f'_*$ it follows
immediately that the homomorphism $(f'_{\rm et})_*$
in \eqref{re3} is surjective. Now Corollary \ref{cor1} says that
$$
{\rm rank}({\mathcal F})\,=\, 1\, ,
$$
where ${\mathcal F}\,\longrightarrow\, X$ is the holomorphic vector bundle
in Proposition \ref{prop2}.

To prove the converse, assume that
\begin{equation}\label{f1}
{\rm rank}({\mathcal F})\,=\, 1\, .
\end{equation}

In view of \eqref{f1}, from Corollary \ref{cor1} it follows that
the homomorphism $(f'_{\rm et})_*$ in \eqref{re3} is surjective. From
the surjectivity of $(f'_{\rm et})_*$ it can be deduced that
the homomorphism of topological fundamental groups
\begin{equation}\label{a3}
f'_*\, :\, \pi^t_1(Y')\, \longrightarrow\, \pi^t_1(\mathcal{X})
\end{equation}
induced by $f'$ in \eqref{e13} is surjective. To see this, first
note that $\pi^t_1(\mathcal{X})$ is residually finite and
$\pi^t_1(Y')$ is finitely generated as they are both surface groups.
Now a result of Peter Scott, \cite[p.~555, Theorem 3.3]{Sc}, says
that for any finitely generated subgroup $H$ of $\pi^t_1(\mathcal{X})$,
and any $t\, \in\, \pi^t_1(\mathcal{X})\setminus H$, there is a finite
index subgroup $$\widetilde{H}\, \subset\, \pi^t_1(\mathcal{X})$$ such that
$$
t\, \notin\, \widetilde{H} \, \supset\, H
$$
(see \cite[p.~2892 Theorem 1.2]{Pat} for an effective version of the theorem
of Scott). Applying this to the image $f'_*(\pi^t_1(Y'))\, \subset\,
\pi^t_1(\mathcal{X})$ we conclude that if $f'_*$ is not surjective then
the image $f'_*(\pi^t_1(Y'))\, \subset\, \pi^t_1(\mathcal{X})$
is contained in a proper subgroup
\begin{equation}\label{t1}
\Gamma\, \subsetneq\, \pi^t_1(\mathcal{X})
\end{equation}
of finite index.

Consider the finite \'etale covering 
$$
\varphi\,\,:\,\, Z\,\, \longrightarrow\, \,{\mathcal X}
$$
given by the subgroup $\Gamma$ in \eqref{t1}. Since
$f'_*(\pi^t_1(Y'))\, \subset\, \Gamma$, 
and a morphism $\beta\,:\, Y\, \longrightarrow\, Z$, such that
$\varphi_0\circ\beta\,=\, f$. But this implies that
$$
(f'_{\rm et})_*(\pi_1(Y'))\, \subset\, \pi_1(Z)
\, \subsetneq\, \pi_1(\mathcal{X})\, ,
$$
where $(f'_{\rm et})_*$ is the homomorphism in \eqref{re3}.
But this contradicts the fact that the homomorphism $(f'_{\rm et})_*$ is surjective.
Therefore, the homomorphism $f'_*$ in \eqref{a3} is surjective.
\end{proof}

\subsection{Pullback of parabolic bundles}

Let $E_*\,=\, (E,\, (\{E_{i,j}\},\, \{\alpha_{i,j}\}))$ be a parabolic vector bundle with
parabolic divisor $D\,=\, \{x_1,\, \cdots, \, x_n\}$. Take a nonconstant holomorphic map
$$
f\, :\, Y\, \longrightarrow\, X
$$
from an irreducible complex projective curve $Y$. Then, using $f$, the parabolic bundle $E_*$
pulls back to a parabolic bundle $f^*E_*$ on $Y$. We will briefly recall the construction
of the parabolic bundle $f^*E_*$.

We first consider the case where $\text{rank}(E)\,=\, 1$. So for each $x_i\, \in\, D$ the 
parabolic weight of $E_*$ is $\alpha_{i,1}\,=\, \alpha_i$. The parabolic divisor for $f^*E_*$ is the 
reduced effective divisor $f^{-1}(D)_{\rm red}$. For $1\, \leq\, i\, \leq\, n$, let 
$$f^{-1}(x_i)\,=\, \{y_{i,1},\, \cdots,\, y_{i,b_i}\}\, \subset\, Y$$ be the inverse 
image, and let $m_{i,j}$ be the multiplicity of $f$ at $y_{i,j}$ for every $1\, \leq\, 
j\, \leq\, b_i$. For any $\lambda\, \in\, \mathbb Q$, let $\lfloor{\lambda}\rfloor$ be 
the integral part of $\lambda$, so $0\, \leq\, \lambda-\lfloor{\lambda}\rfloor \, <\, 1$.

The holomorphic line bundle on $Y$ underlying the parabolic line bundle $f^*E_*$ is
$$
F\, :=\, (f^*E)\otimes{\mathcal O}_Y(\sum_{i=1}^n\sum_{j=1}^{b_i} \lfloor{m_{i,j}\alpha_i}\rfloor
\cdot y_{i,j})\, ,
$$
and the parabolic weight of $F_{y_{i,j}}$ is $m_{i,j}\alpha_i -
\lfloor{m_{i,j}\alpha_i}\rfloor$. Note that
\begin{equation}\label{d1}
\text{par-deg}(f^*E_*)\,=\, \text{degree}(f)\cdot \text{par-deg}(E_*)\, .
\end{equation}

Any parabolic vector bundle $E_*$ can locally be expressed as a direct sum of parabolic 
line bundles. In other words, $X$ can be covered by Zariski open subsets $U_1,\, 
\cdots,\, U_m$ such that $E_*\big\vert_{U_j}$ is a direct sum of parabolic line bundles 
on $U_j$ for all $1\, \leq\, j\, \leq\, m$. We have described above the pullback of 
parabolic line bundles. The pullback of a direct sum of parabolic line bundles is the 
direct sum of the pulled back parabolic line bundles. Using the decomposition of 
$E_*\big\vert_{U_j}$ into a direct sum of parabolic line bundles we now have a 
description of the parabolic pullback $f^*E_*$. From \eqref{d1} it follows that
\begin{equation}\label{d2}
\text{par-deg}(f^*E_*)\,=\, \text{degree}(f)\cdot \text{par-deg}(E_*)
\end{equation}
for any parabolic vector bundle $E_*$.

For each point $x\, \in\, D$ fix an integer $N_x\, \geq\, 1$.

\begin{theorem}\label{thm2}
Take any stable parabolic vector bundle $E_*$ on $X$ with parabolic structure over $D$ 
such that all the parabolic weights of $E_*$ at each $x\, \in\, D$ are integral multiples 
of $\frac{1}{N_x}$. If the rank of the holomorphic vector bundle ${\mathcal 
F}\,\longrightarrow\, X$ in Proposition \ref{prop2} is one, then the parabolic vector 
bundle $f^*E_*$ on $Y$ is also stable.
\end{theorem}

\begin{proof}
Assume that
\begin{equation}\label{d3}
{\rm rank}(\mathcal{F})\,=\, 1\, .
\end{equation}

Let $E_*$ be a stable parabolic vector bundle of rank $r$ on $X$ with parabolic structure 
over $D$ such that all the parabolic weights of $E_*$ at each $x\, \in\, D$ are integral 
multiples of $\frac{1}{N_x}$. Consider the parabolic principal $\text{PGL}(r, {\mathbb 
C})$--bundle $\mathbb{P}(E_*)$ defined by $E_*$; see \cite{BBN} for parabolic principal 
bundles. Since $E_*$ is stable, we know that the parabolic principal $\text{PGL}(r, 
{\mathbb C})$--bundle $\mathbb{P}(E_*)$ is given by an irreducible homomorphism
\begin{equation}\label{r1}
\rho\,\, :\,\, \pi^t_1(\mathcal{X})\, \longrightarrow\, \text{PU}(r)
\end{equation}
\cite{MS}, \cite{Biq}. Let $\mathbb{P}(f^*E_*)$ denote the parabolic principal 
$\text{PGL}(r, {\mathbb C})$--bundle on $Y$ defined by the parabolic vector bundle 
$f^*E_*$. Since $\mathbb{P}(E_*)$ is given by the homomorphism $\rho$ in \eqref{r1}, we 
conclude that the parabolic principal $\text{PGL}(r, {\mathbb C})$--bundle 
$\mathbb{P}(f^*E_*)$ is given by the homomorphism
\begin{equation}\label{r2}
\rho\circ f'_* \,\, :\,\, \pi^t_1(Y')\, \longrightarrow\, \text{PU}(r)\, ,
\end{equation}
where $f'_*$ is the homomorphism in \eqref{r3}.

{}From \eqref{d3} and Proposition \ref{prop3} we know that the homomorphism $f'_*$ in 
\eqref{r2} is surjective. Therefore, from the property of the homomorphism $\rho$ in 
\eqref{r1} that it is irreducible we conclude that the homomorphism $\rho\circ f'_*$ in 
\eqref{r2} is also irreducible. Since the parabolic principal $\text{PGL}(r, {\mathbb 
C})$--bundle $\mathbb{P}(f^*E_*)$ is given by the irreducible projective unitary 
representation $\rho\circ f'_*$ in \eqref{r2}, we now conclude that the parabolic vector 
bundle $f^*E_*$ is stable \cite{MS}, \cite{Biq}.
\end{proof}

The following lemma is a converse of Theorem \ref{thm2}.

\begin{lemma}\label{lem4}
Assume that the rank of the holomorphic vector bundle ${\mathcal
F}\,\longrightarrow\, X$ in Proposition \ref{prop2} is at least two.
Then there is a stable parabolic vector bundle $E_*$ on $X$ with parabolic structure over $D$
such that
\begin{enumerate}
\item all the parabolic weights of $E_*$ at each $x\, \in\, D$ are integral multiples
of $\frac{1}{N_x}$, and

\item the parabolic vector bundle $f^*E_*$ on $Y$ is not stable.
\end{enumerate}
\end{lemma}

\begin{proof}
Since ${\rm rank}({\mathcal F})\, >\, 1$, from Proposition \ref{prop3}
we know that the homomorphism of topological fundamental groups
$$
f'_* \, :\, \pi^t_1(Y')\, \longrightarrow\, \pi^t_1(\mathcal{X})
$$
(see \eqref{r3}) induced by $f'$ in \eqref{e13} is not surjective. Fix an
irreducible representation
$$
\rho\, :\, \pi^t_1(\mathcal{X})\, \longrightarrow\, {\rm U}(r)\, ,
$$
for some $r\, \geq\, 2$, such that the composition of homomorphisms
$$
\rho\circ f'_* \, :\, \pi^t_1(Y')\, \longrightarrow\, {\rm U}(r)
$$
is not irreducible; such a $\rho$ exists because $f'_*$ is not surjective.

Let $E_*$ be the parabolic vector bundle of rank $r$ on $X$, with parabolic structure over $D$,
given by $\rho$. We note that
\begin{enumerate}
\item all the parabolic weights of $E_*$ at each $x\, \in\, D$ are integral multiples
of $\frac{1}{N_x}$, and

\item the parabolic vector bundle $E_*$ is stable, because $\rho$ is irreducible \cite{MS}.
\end{enumerate}
Since $\rho\circ f'_* \, :\, \pi^t_1(Y')\, \longrightarrow\, {\rm U}(r)$ is not irreducible,
it follows that the parabolic vector bundle $f^*E_*$ on $Y$ is not stable.
\end{proof}

\section{Algebraically closed fields of characteristic zero}

Now let $k$ be any algebraically closed fields of characteristic zero. As in Section \ref{se3},
$X$ is an irreducible smooth projective curve defined over $k$, and
$$
D\, =\, \{x_1,\, \cdots ,\, x_n\}\, \subset\, X
$$
is a finite subset.

For each point $x\, \in\, D$ fix an integer $N_x\, \geq\, 1$.
Consider the holomorphic vector bundle ${\mathcal F}\,\longrightarrow\, X$ 
in Proposition \ref{prop2}.

\begin{theorem}\label{thm3}
Let
$$
f\, :\, Y\, \longrightarrow\, X
$$
be a nonconstant map from an irreducible smooth projective curve $Y$. Take any stable 
parabolic vector bundle $E_*$ on $X$ with parabolic structure over $D$ such that all the 
parabolic weights of $E_*$ at each point $x\, \in\, D$ are integral multiples of 
$\frac{1}{N_x}$. If the rank of the vector bundle ${\mathcal F}\,\longrightarrow\, X$ in 
Proposition \ref{prop2} is one, then the parabolic vector bundle $f^*E_*$ on $Y$ is also 
stable.
\end{theorem}

\begin{proof}
Assume that
\begin{equation}\label{l1}
\text{rank}(\mathcal{F})\,=\, 1\, .
\end{equation}

Take any stable parabolic vector bundle $E_*$ on $X$ with parabolic structure over $D$ 
such that all the parabolic weights of $E_*$ at each $x\, \in\, D$ are integral multiples 
of $\frac{1}{N_x}$. We need to show that the parabolic vector bundle $f^*E_*$ on $Y$ is 
stable.

Let $k_0 \,\subset\, k$ be an algebraically closed field of characteristic of finite 
transcendence degree over $\mathbb Q$ such that $X$, $Y$, $D$, $f$ and $E_*$ are defined 
over $k_0$. Fix an embedding of $k_0$ in $\mathbb C$. Let
\begin{equation}\label{b1}
X_{\mathbb C}\, :=\, X\times_{k_0} {\mathbb C},\ Y_{\mathbb C}\, :=\, Y\times_{k_0} 
{\mathbb C},\ f_{\mathbb C}\, :=\, f\times_{k_0} {\mathbb C}\,\, \text{ and }\,\, 
E^{\mathbb C}_*\, :=\, E_*\otimes_{k_0} {\mathbb C}
\end{equation}
be the base changes to $\mathbb C$ of $X$, $Y$, $f$ and $E_*$ respectively. Similarly, let
\begin{equation}\label{b2}
\mathcal{X}_{\mathbb C}\, :=\, \mathcal{X}\times_{k_0} {\mathbb C},\ 
Y'_{\mathbb C}\, :=\, Y'\times_{k_0} {\mathbb C}\ \text{ and }\ f'_{\mathbb C}\,
:=\, f'\times_{k_0} {\mathbb C}
\end{equation}
where $Y'$ and $f'$ are as in \eqref{e13}, be the base changes to $\mathbb C$ of 
$\mathcal X$, $Y'$ and $f'$ respectively.

We need the following lemma.

\begin{lemma}\label{lem3}
The parabolic vector bundle $E^{\mathbb C}_*\, =\, E_*\otimes_{k_0} {\mathbb C}$
in \eqref{b1} is stable.
\end{lemma}

\begin{proof}[{Proof of Lemma \ref{lem3}}]
An equivariant vector bundle is equivariantly semistable if the underlying vector bundle 
is semistable, because the Harder--Narasimhan filtration of an equivariant bundle is 
preserved by the action of the group. It is known that the property of semistability of a 
vector bundle is preserved under field extensions (see \cite[p.~18, Corollary 
1.3.8]{HL}). Now using the correspondence between the parabolic bundles and the 
equivariant bundles we conclude that for a semistable parabolic bundle $V_*$ on $X$ the 
parabolic bundle $V_*\times_{k_0} {\mathbb C}$ on $X_{\mathbb C}$ is also semistable. 
Therefore, the given condition that the parabolic bundle $E_*$ is semistable implies that 
the parabolic bundle $E^{\mathbb C}_*$ is also semistable. Since the unique maximal 
polystable parabolic subbundle of the semistable parabolic bundle $E^{\mathbb C}_*$ (it 
is also known as the socle of $E^{\mathbb C}_*$ (see \cite[p.~23, Lemma 1.5.5]{HL}) is 
defined over $k_0$, and $E_*$ is polystable, we conclude that the parabolic bundle 
$E^{\mathbb C}_*$ is polystable (see \cite[p.~24, Corollary 1.5.11]{HL}.

For a parabolic vector bundle $F_*$, the sheaf of quasiparabolic structure preserving 
endomorphisms of the underlying vector bundle $F$ will be denoted by $\text{End}_P(F_*)$. 
A polystable parabolic vector bundle $F_*$ is stable if and only if the space of global 
sections of $\text{End}_P(F_*)$ is the base field. Since
$$
\text{End}_P(E^{\mathbb C}_*)\,=\, \text{End}_P(E_*) \otimes_{k_0} {\mathbb C}\, ,
$$
and $E_*$ is stable, we have
$$H^0(X_{\mathbb C},\, \text{End}_P(E^{\mathbb C}_*))\,=\, H^0(X,\, \text{End}_P(E_*)) \otimes_{k_0} {\mathbb C}
\,=\, \mathbb C\, .$$
This implies that the polystable parabolic bundle $E^{\mathbb C}_*$ is stable.
\end{proof}

Continuing with the proof of Theorem \ref{thm3}, from Corollary \ref{cor1} and \eqref{l1} 
we know that the homomorphism of \'etale fundamental groups $(f'_{\rm et})_*$ in 
\eqref{re3} is surjective. This implies that the homomorphism of \'etale fundamental 
groups
$$
(f'_{\mathbb{C}, {\rm et}})_*\, :\, \pi_1(Y'_{\mathbb C})\, \longrightarrow\,
\pi_1(\mathcal{X}_{\mathbb C})
$$
induced by $f'_{\mathbb C}$ in \eqref{b2} is surjective; both $Y'_{\mathbb C}$ and 
$\mathcal{X}$ are defined in \eqref{b2}. Hence from Theorem \ref{thm2} and Lemma 
\ref{lem3} we conclude that the parabolic vector bundle $f^*_{\mathbb C}E^{\mathbb C}_*$ 
on $Y_{\mathbb C}$ is stable, where $f_{\mathbb C}$ is the map in \eqref{b1}. This 
implies that the parabolic vector bundle $f^*E$ is stable.
\end{proof}

\begin{remark}
We note that the main result of \cite{GR-1} implies that if the base field is algebraically 
closed of characteristic 0 then the conclusion to Theorem \ref{thm3} holds
under the strict condition that for 
every $x\,\in\, X$ the number $N_x$ is coprime to the ramification indices of $f$ at points 
above $x$. Theorem \ref{thm3} is more general than this as illustrated by Example \ref{exm}.
\end{remark}

\begin{example}\label{exm}
Let $X\,=\,Y\,=\,\mathbb{P}^1_{\mathbb C}$ and $f\,:\,Y\,\longrightarrow\, X$
be the cyclic covering of degree 
6 ramified at $0$ and $\infty$. Let $N_0\,=\,2$ and $N_{\infty}\,=\,3$; then the map 
$\pi_1(Y\setminus \{0,\,\infty\})\,\longrightarrow\, \pi_1(\mathcal{X})$ is surjective, hence
the rank of $\mathcal{F}$ is one so Theorem \ref{thm3} applies. Though the hypothesis of
\cite[Theorem 5.1]{GR-1} does not hold in this example.
\end{example}

\section*{Acknowledgements}

We are very grateful to the referee for the question on the proof of Proposition \ref{prop3}.
We thank Mahan Mj for pointing out \cite{Sc} and \cite{Pat}. 
The first-named author thanks Universit\'e C\^ote d'Azur for hospitality while parts of the
work were carried out. He is partially supported by a J. C. Bose Fellowship.



\begin{thebibliography}{ZZZZ}

\bibitem[BBN]{BBN} V. Balaji, I. Biswas and D. S. Nagaraj, Principal bundles over 
projective manifolds with parabolic structure over a divisor, {\it Tohoku Math. J.}
{\bf 53} (2001), 337--367.

\bibitem[Bh]{Bh} U. N. Bhosle, Parabolic sheaves on higher-dimensional varieties,
{\it Math. Ann.} {\bf 293} (1992), 177--192.

\bibitem[Biq]{Biq} O. Biquard, Fibr\'es paraboliques stables et connexions singuli\`eres
plates, {\it Bull. Soc. Math. Fr.} {\bf 119} (1991), 231--257.

\bibitem[Bis1]{Bi1} I. Biswas, Parabolic bundles as orbifold bundles, {\it Duke Math. J.}
{\bf 88} (1997), 305--325.

\bibitem[Bis2]{Bi2} I. Biswas, A cohomological criterion for semistable parabolic vector 
bundles on a curve, {\it Com. Ren. Math. Acad. Sci. Paris} {\bf 345} (2007), 325--328.

\bibitem[BP]{BP} I. Biswas and A. J. Parameswaran, Ramified covering maps and stability of 
pulled back bundles, {\it Int. Math. Res. Not.}, https://doi.org/10.1093/imrn/rnab062.

\bibitem[BKP]{GR-1} I. Biswas, M.Kumar and A. J. Parameswaran, Genuinely ramified maps and 
stability of pulled-back parabolic bundles, {\it Indag. Math.},
https://doi.org/10.1016/j.indag.2022.04.003.

\bibitem[Bo1]{Bo1} N. Borne, Fibr\'es paraboliques et champ des racines, {\em Int. Math.
Res. Not. IMRN}, {\bf 16}, Art. ID rnm049, 38, (2007).

\bibitem[Bo2]{Bo2} N. Borne, Sur les repr\'esentations du groupe fondamental d'une
vari\'et\'e priv\'ee d'un diviseur \`a croisements normaux simples, {\em Indiana Univ.
Math. Jour.} {\bf 58} (2009), 137--180.

\bibitem[Ha]{Ha} R. Hartshorne, {\it Algebraic geometry}, Graduate Texts in Mathematics,
No. 52. Springer-Verlag, New York-Heidelberg, 1977.

\bibitem[HL]{HL} D. Huybrechts and M. Lehn, {\it The geometry of moduli spaces of 
sheaves}, Aspects of Mathematics, E31, Friedr. Vieweg~\&~Sohn, Braunschweig, 1997.

\bibitem[IIS]{IIS} M.-a. Inaba, K. Iwasaki and M.-H. Saito, Moduli of stable parabolic
connections, Riemann-Hilbert correspondence and geometry of Painlev\'e equation of type VI. I,
{\it Publ. Res. Inst. Math. Sci.} {\bf 42} (2006), 987--1089.

\bibitem[In]{In} M.-a. Inaba, Moduli of parabolic connections on curves and the Riemann-Hilbert
correspondence, {\it J. Algebraic Geom.} {\bf 22} (2013), 407--480.

\bibitem[MY]{MY} M. Maruyama and K. Yokogawa, Moduli of parabolic stable sheaves, Math.
Ann. 293 (1992), 77--99.

\bibitem[MS]{MS} V. B. Mehta and C. S. Seshadri, Moduli of vector bundles on curves with
parabolic structures, \textit{Math. Ann.} \textbf{248} (1980), 205--239.

\bibitem[Par]{Pa} A. J. Parameswaran, Parabolic coverings I: the case of curves,
{\it J. Ramanujan Math. Soc.} {\bf 25} (2010), 233--251.

\bibitem[Pat]{Pat} P. Patel, On a theorem of Peter Scott, {\it Proc. Amer. Math. Soc.}
{\bf 142} (2014), 2891--2906.

\bibitem[Sc]{Sc} P. Scott, Subgroups of surface groups are almost geometric,
{\it Journal London Math. Soc.} {\bf 2} (1978), 555--565.

\bibitem[Yo]{Yo} K. Yokogawa, Infinitesimal deformations of parabolic Higgs sheaves,
{\it Internet. J. Math.} {\bf 6} (1995), 125--148. 

\end{thebibliography}
\end{document}